\documentclass[12pt,a4paper]{article}

\usepackage{a4wide,amssymb,amsmath,amsthm,xspace,epsfig, amsfonts}

\usepackage[left=1.5cm, right=1.5cm, top=1.5cm]{geometry}

\begin{document}

\newcommand{\N}{\mathbb{N}}
\newcommand{\R}{\mathbb{R}}
\newcommand{\Z}{\mathbb{Z}}
\newcommand{\Q}{\mathbb{Q}}
\newcommand{\C}{\mathbb{C}}
\newcommand{\PP}{\mathbb{P}}

\newcommand{\LL}{\Bbb L}
\newcommand{\OO}{\mathcal{O}}

\newcommand{\esp}{\vskip .3cm \noindent}
\mathchardef\flat="115B

\newcommand{\lev}{\text{\rm Lev}}

\def\ut#1{$\underline{\text{#1}}$}
\def\CC#1{${\cal C}^{#1}$}
\def\h#1{\hat #1}
\def\t#1{\tilde #1}
\def\wt#1{\widetilde{#1}}
\def\wh#1{\widehat{#1}}
\def\wb#1{\overline{#1}}

\def\restrict#1{\bigr|_{#1}}

\def\ufin#1#2{\mathsf{U_{fin}}\bigl({#1},{#2}\bigr)}
\def\sfin#1#2{\mathsf{S_{fin}}\bigl({#1},{#2}\bigr)}
\def\s1#1#2{\mathsf{S_{1}}\bigl({#1},{#2}\bigr)}

\def\ch#1#2{\left(\begin{array}{c}#1 \\ #2 \end{array}\right)}

\newtheorem{lemma}{Lemma}[section]

\newtheorem{thm}[lemma]{Theorem}

\newtheorem{defi}[lemma]{Definition}
\newtheorem{conj}[lemma]{Conjecture}
\newtheorem{cor}[lemma]{Corollary}
\newtheorem{prop}[lemma]{Proposition}
\newtheorem{prob}{Problem}
\newtheorem{q}[lemma]{Question}
\newtheorem*{rem}{Remark}
\newtheorem{examples}[lemma]{Examples}
\newtheorem{example}[lemma]{Example}

\title{Three small results on normal first countable linearly H-closed spaces}
\date{\empty}
\author{Mathieu Baillif}
\maketitle

\abstract{\footnotesize{We use topological consequences of {\bf PFA}, {\bf MA$_{\omega_1}$(S)[S]} and {\bf PFA(S)[S]} proved by other authors
to show that normal first countable linearly H-closed spaces with various additionals properties are compact in these models.}}

\section{Statements}

In this small note, we prove three modest results about the following problem we raised in \cite[Question 2.12]{meszigues-od-sel}, 
using topological consequences of the axioms {\bf MA$_{\omega_1}$(S)[S], PFA(S)[S]} and {\bf PFA} due to other authors.

\begin{prob}\label{prob:1}
   Is there in {\bf ZFC} a non-compact normal linearly H-closed space which is first countable (or equivalently has $G_\delta$-points)~?
\end{prob}

By `space' we mean `topological Hausdorff space', hence `regular' and `normal' imply `Hausdorff'. 
A {\em cover} of a space always means a cover by open sets, and a cover is a {\em chain cover} if it is linearly ordered by the inclusion relation.
A space is {\em linearly H-closed} provided any chain cover has a dense member. There are various equivalent definitions,
see for instance Lemmas 2.2--2.3 in \cite{meszigues-od-sel} or Theorems 2.5, 2.12--2.14 in \cite{AlasJunqueiraWilson:2019}.
A space is {\em H-closed} iff any cover has a finite subfamily with a dense union. Regular H-closed spaces are compact.

We strongly suspect that solving Problem \ref{prob:1} is not listed very high on mankind's priorities list
but 
hope that our results are of some interest for researchers in the field. 
To help enhancing their curiosity (and for context), we start by listing known relevant results.
Recall that a space is {\em feebly compact} iff every locally finite family of non-empty open sets is finite or equivalently
(for Hausdorff spaces)
iff any countable cover has a a finite subfamily with a dense union. A linearly H-closed space is feebly compact. 
A feebly compact space 
is {\em pseudocompact} (that is, any real valued function is bounded)  and the converse holds for Tychonoff spaces.
(We do not include regularity in the definitions of pseudocompactness and feeble compactness.)
Recall also that a pseudocompact normal space is 
countably compact. 

In the list below,
we tried to give the name of the person who first constructed the example 
(in general in contexts unrelated to linearly H-closed spaces) or proved the result. 
Each fact given without a reference is explained in \cite[Example 2.9]{meszigues-od-sel}, where links to
the original sources can be found.
We assume that the reader knows what the axioms 
{\bf CH}, {\bf MA + $\neg$CH}, {\bf PFA}, $\diamondsuit$, $\mathfrak{p}=1$ mean.
For undefined terms see the given references. 
Fact \ref{fact:Gdelta} explains the equivalence
in the statement of Problem \ref{prob:1}.

\vskip .3cm\noindent
\underline{A list of known facts} 
\begin{enumerate}
\renewcommand\theenumi{\Alph{enumi}}
\item \label{fact:Gdelta}
       A regular feebly compact space with $G_\delta$ points is first countable (Porter-Woods, \cite[Prop. 2.2]{PorterWoods:1984}).
\item  A weakly linearly Lindel\"of (in particular weakly Lindel\"of or ccc) feebly compact space is linearly H-closed
       (Alas, Junqueira and Wilson \cite[Thm. 2.13]{AlasJunqueiraWilson:2019}).
\item \label{fact:justmissed}
      There are non-compact linearly H-closed spaces which are moreover
      perfect, Tychonoff and first countable (for instance the famous space $\Psi$
       due to Isbell and Mr\'owka), 
       or Frechet-Urysohn and collectionwise normal (for instance a $\Sigma$-product of $2^{\omega_1}$).
\item Under {\bf CH}, a normal first countable linearly H-closed space has cardinality $\le\aleph_1$ and is weakly Lindel\"of (A. Bella 
      \cite[Thm 4.4 \& Cor. 4.6]{Bella:2017}).
\item \label{fact:p=omega_1}
          If $\mathfrak{p} = \omega_1$ (in particular under {\bf CH}), there is a non-compact 
          normal separable first countable locally compact locally countable linearly H-closed space, namely a space of the type $\gamma\N$
          due originally to Franklin and Rajagopalan and studied in detail by other authors, in particular by Nyikos. 
\item \label{fact:diamond} Under $\diamondsuit$, there is a non-compact perfectly normal first countable hereditarily separable
          linearly H-closed space, for instance an Ostaszewski space, 
          and even a manifold (M.E. Rudin, see Nyikos' exposition in  \cite[Ex. 3.14]{Nyikos:1984}).
          If one adds Cohen reals to the model, the manifold keeps these properties \cite{BaloghGruenhage:2005}.
\item Under {\bf MA + $\neg$CH}, a perfectly normal linearly H-closed space is compact. This result is
            also compatible with {\bf CH}. (Weiss and Eisworth, see
            \cite[Lemma 2.1]{meszigues-od-sel}.) 
\item \label{fact:PFA}
          Under {\bf PFA}, a normal linearly H-closed space which is either locally separable and countably tight or locally ccc, locally compact and 
          first countable is compact.
          \cite[Thm 2.13]{meszigues-od-sel}.
\item \label{fact:mononoke} A monotonically normal linearly H-closed space is compact (Alas, Junqueira and Wilson \cite[Thm. 2.17]{AlasJunqueiraWilson:2019}).
\end{enumerate}

Our new results are the following. 
They fit between Facts \ref{fact:PFA} and \ref{fact:mononoke}, so to say.
Recall that the {\em spread} of a space is the supremum of the cardinalities of its discrete subspaces.
Hereditarily separable spaces (such as those in Fact \ref{fact:diamond}) are of countable spread.

\begin{thm}
   \label{thm:PFAspread}
   {\bf PFA} implies that if $X$ is a normal linearly H-closed space which is locally of countable spread, then
   $X$ is compact and first countable.
\end{thm}

\begin{thm}\label{thm:FCLC}
   In a particular model of {\bf MA$_{\omega_1}$(S)[S]} and in
   any model of {\bf PFA(S)[S]},   
   a locally compact, hereditarily normal, linearly H-closed space with $G_\delta$ points is compact.
\end{thm}

\begin{thm}
  \label{thm:PFASS}
  In any model of {\bf PFA(S)[S]},
  a linearly H-closed hereditarily normal space such that each point is a $G_\delta$ and has
  an open Lindel\"of neighborhood is compact.
\end{thm}

{\bf MA$_{\omega_1}$(S)[S]} and {\bf PFA(S)[S]} are formally not forcing axioms but rather a powerful method for obtaining models
starting from a coherent Suslin tree $S$, using iterated forcing to obtain weaker versions (called {\bf MA$_{\omega_1}$(S)} and {\bf PFA(S)})
of 
{\bf MA$_{\omega_1}$} and {\bf PFA} which preserve (the Suslinity of) $S$,
and then forcing with $S$.
The precise definitions shall not concern us here as we use only topological implications which hold
in these models. 
The known proof of the consistency of {\bf PFA(S)} needs inaccessible cardinals, while that of 
{\bf MA$_{\omega_1}$(S)} does not. Hence when we write `in a model of {\bf MA$_{\omega_1}$(S)[S]}'
it is implied that the model is obtained without inaccessibles.
For details, see F.D. Tall article \cite{Tall:PFAforthemasses} and references therein.
We recall that {\bf PFA} $\Longrightarrow$ {\bf MA + $\neg$CH}
and 
{\bf PFA(S)[S]} $\Longrightarrow$ {\bf MA$_{\omega_1}$(S)[S]}.

We note that {\bf PFA(S)[S]} implies $\mathfrak{p}=\omega_1$ (see the introduction of \cite{Tall:PFAforthemasses}), 
hence there is a model of set theory where Theorems \ref{thm:FCLC}--\ref{thm:PFASS} and Fact \ref{fact:p=omega_1} hold,
that is, there are first countable linearly H-closed locally compact spaces which are normal, but none which is hereditarily normal.

\section{Proofs}

Our proofs are short and quite similar to patchworks, as we mostly blend together results found elsewhere.
Since reading a sequence of references can be somewhat dull, we provide proofs (if short enough)
of the results that do not appear as an explicit lemma or theorem
somewhere else, or for which we have small variants in the argument. Let us start with the list.
First,
the following easy fact will be used several times.
\begin{lemma}[{\cite[Lemma 2.3]{meszigues-od-sel} or \cite[Thm 2.1]{AlasJunqueiraWilson:2019}}]
   \label{lemma:openLHC}
   If $X$ is linearly H-closed and $U\subset X$ is open, then $\wb{U}$ is linearly H-closed.
\end{lemma}

We denote the following result by {\bf HL} following \cite{Tall:PFAforthemasses}.
Recall that an {\em S-space} is a regular hereditarily separable non-Lindel\"of space,
while an {\em L-space} is a regular hereditarily Lindel\"of  non-separable space.

\begin{lemma}[{\cite{Szentmiklossy:1980}, \cite[Thm 2.1]{AbrahamTodorcevic:1984}, \cite[Lemma 11]{Tall:PFAArhan}}]
   \label{lemma:MAorMA}
   {\bf (MA$_{\omega_1}$(S)[S]} or {\bf MA + $\neg$CH)}\
   \\
   {\bf HL} 
   A first countable hereditarily Lindel\"of regular space is hereditarily separable, that is,
   first countable L-spaces do not exist.
\end{lemma}

We shall also use the basic result below which is part of the folklore and is cited for instance in \cite{EisworthNyikosShelah:2003}.
A detailed proof can be found online in Dan Ma's topology blog\footnote{Available at 
{\em https://dantopology.wordpress.com/2018/10/15/a-little-corner-in-the-world-of-set-theoretic-topology/} when this note was written.}.

\begin{lemma}\label{lemma:spreadL} 
   If a regular space of countable spread is not hereditarily separable, it contains an L-space, 
   and if it is not hereditarily Lindelof, it contains an S-space.
\end{lemma}

The next lemma is due to Todor\v cevi\'c. 

\begin{lemma}[{Todor\v cevi\'c, \cite[Theorem 8.11]{TodorcevicBook}}]
   \label{lemma:JuPFA} \ \\
   {\bf (PFA)}
   A space of countable spread has $G_\delta$ points.
\end{lemma}
It then follows that:
\begin{lemma}[{\bf PFA}]
   \label{lemma:regcount}
   If $X$ is a regular linearly H-closed space which is locally of countable spread, then
   $X$ is first countable and locally hereditarily separable.
\end{lemma}

\proof
   Let $U_x\ni x$ be an open neighborhood of countable spread.
   Since $X$ is regular, we may choose an open $V_x$ such that
   $x\in V_x\subset\wb{V_x}\subset U_x$. 
   Then $\wb{V_x}$ is linearly H-closed by Lemma \ref{lemma:openLHC}.
   Then by Lemma \ref{lemma:JuPFA} $\wb{V_x}$ has $G_\delta$ points under {\bf PFA} and hence is first countable by Fact \ref{fact:Gdelta}.
   If $\wb{V_x}$ is not hereditarily separable by Lemma \ref{lemma:spreadL} it contains an L-space, but 
   since {\bf PFA} implies {\bf MA + $\neg$CH}, {\bf HL} (Lemma \ref{lemma:MAorMA}) shows that it is impossible.
   Hence $\wb{V_x}$ is hereditarily separable.
\endproof

Theorem \ref{thm:PFAspread} follows almost immediately.
\proof[Proof of Theorem \ref{thm:PFAspread}]
   Let $X$ be normal, linearly H-closed and locally of countable spread. By Lemmas \ref{lemma:openLHC} and \ref{lemma:regcount},
   any open subset of countable spread is separable and $X$ is first countable.
   Hence $X$ is locally separable, 
   which implies the result by Fact \ref{fact:PFA}. 
\endproof

We now turn our attention to the other two theorems.
We need the following easy lemma.

\begin{lemma}\label{lemma:loccomp}
  If $X$ is locally compact, countably tight, linearly H-closed and non-compact, then
  $X$ contains a $\sigma$-compact open set $U$ such that $\wb{U}$ is linearly
  H-closed and non-compact.
\end{lemma}
Recall that a space is {\em countably tight} iff for each subset $E$ and each point $x$ in its closure, there
is a countable $A\subset E$ whose closure contains $x$.
\proof
 First, notice that a linearly H-closed Lindel\"of space is H-closed and hence compact if regular.
 Hence  by Lemma \ref{lemma:openLHC} the closure of an open set in $X$ is either compact or non-Lindel\"of.
 We build open sets $U_\alpha$ with compact closure (indexed by countable ordinals) such that $\wb{U_\alpha}\subset U_\beta$ whenever $\alpha<\beta$
 by induction, starting with some $U_0$.
 Given $U_\beta$ for each $\beta<\alpha$, either $\wb{\cup_{\beta<\alpha} U_\beta}$ is non-Lindel\"of and
 we are over, or it is compact. In the latter case,
 since $X$ is non-compact, we may choose a point $x_\alpha\not\in\wb{\cup_{\beta<\alpha}U_\beta}$.
 We then cover $\{x_\alpha\}\cup\wb{\cup_{\beta<\alpha} U_\beta}$ with finitely many open sets with compact closure whose
 union defines $U_\alpha$. Notice that $\wb{U_\alpha}\subsetneq U_\beta$ when $\alpha<\beta$.
 If this goes on until $\omega_1$, then $W=\cup_{\alpha<\omega_1} U_\alpha$ is clopen.
 Indeed, openness is immediate, and by countable tightness any point $x\in\wb{W}$ is in the closure of a countable subset of $W$ which must be 
 contained in some $U_\alpha$, so $x\in\wb{U_\alpha}\subset U_{\alpha+1}\subset W$.
 By Lemma \ref{lemma:openLHC}, $W$ is
 linearly H-closed, which is impossible since none of the $U_\alpha$ is dense.
 Hence the process must stop before $\omega_1$, that is,
 $\wb{\cup_{\beta<\alpha} U_\beta}$ is non-Lindel\"of for some $\alpha<\omega_1$.
\endproof

The next lemma also holds in {\bf ZFC} and is due to Nyikos.
\begin{lemma}[{Nyikos \cite[Lemma 1.2]{Nyikos:2003}}]
  \label{lemma:Nyikos+}
  Let $X$ be a hereditarily $\aleph_1$-scwH space, and $E\subset U\subset X$ be such that $U$ is open, $E$ is $\omega_1$-compact and dense in $U$.
  Then $\wb{U}-U$ has countable spread.
\end{lemma}
Recall that {\em scwH} is a shorthand for {\em strongly collectionwise Hausdorff}. A space is 
$\aleph_1$-scwH iff any closed discrete subset of size $\le\aleph_1$ can be expanded to a discrete family of open sets.
A space is {\em $\omega_1$-compact} (or has {\em countable extent}) iff any closed discrete subspace is at most countable.
Countably compact and Lindel\"of spaces are $\omega_1$-compact.
\proof
   We give a very small variant of Nyikos' proof.
   Let $D$ be a discrete subset of $\wb{U}-U$. Then $D$ is closed discrete in the space $W=D\cup U$, because the closure of $D$ does not intersect $U$.
   Since $W$ is $\aleph_1$-scwH, we may let $\{V_d\,:\,d\in D\}$ be a discrete-in-$W$ collection of open subsets,
   then $\{V_d\cap E \,:\,d\in D\}$ is a discrete-in-$E$ collection of non-empty (by density of $E$) open sets.
   Since $E$ is $\omega_1$-compact, this collection is countable, and so is $D$.
\endproof

We also need some topological implications valid in (some) models of {\bf MA$_{\omega_1}$(S)[S]}.
The first one was originally proved by Szentmikl\'ossy under {\bf MA+$\neg$CH}.
The fact that it holds under {\bf PFA(S)[S]} is due to Todor\v cevi\'c, and in a model of {\bf MA$_{\omega_1}$(S)[S]}
by Larson and Tall, see \cite{LarsonTall:2010b}, \cite[Thm 4.1]{Tall:PFAforthemasses} and references therein.
We provide a short argument  that it follows from $\mathbf{\Sigma^-}$ (defined below, the terminology
and the proof are taken from \cite{LarsonTall:2010b} ).

\begin{thm}[{In a model of {\bf MA$_{\omega_1}$(S)[S]} and in every model of {\bf PFA(S)[S]}}]\label{lemma:MALC}
   A locally compact space of countable spread is hereditarily Lindel\"of. 
\end{thm}
\proof
   It is shown in \cite[Thm 4.1]{Tall:PFAforthemasses} that the following holds in a model of {\bf MA$_{\omega_1}$(S)[S]}
   and in every model of {\bf PFA(S)[S]}:
   \vskip .3cm\noindent
   $\mathbf{\Sigma^-}$: In a compact countably tight space, locally countable subspaces of size $\aleph_1$ are $\sigma$-discrete.
   \vskip .3cm\noindent
   If $X$ is locally compact and of countable spread, its one point compactification $X^*$ is countably tight (and has countable spread). Indeed,
   by a well known fact if $X^*$ is not countably tight, it contains a perfect preimage of $\omega_1$ (see e.g. \cite[Lemma 4]{LarsonTall:2010b}),
   that is, a space $Y$ with a closed onto map $p:Y \to\omega_1$ such that preimages of points are compact. 
   But it is then easy to find an uncountable discrete subspace in $Y$.
   If $X^*$ is not hereditarily Lindel\"of, by Lemma \ref{lemma:spreadL} it contains an S-space and hence (by classical results, 
   see e.g. \cite[Cor. 3.2]{Roitman:1984})
   a subspace $Z=\{x_\alpha\,:\,\alpha\in\omega_1\}$ such that $\{x_\alpha\,:\,\alpha<\gamma\}$ is open in $Z$ for each $\gamma$.
   But then $\mathbf{\Sigma^-}$ implies that $X^*$ contains an uncountable discrete subset, a contradiction.
\endproof

We will also use the following result.
The notation {\bf CW} comes from \cite{Tall:PFAforthemasses} again.
\begin{thm}[{\cite[Cor. 14]{LarsonTall:2010}}]
  \label{thm:N-CWH}
     \ \\
    {\bf CW} After forcing with a Suslin tree, any first countable normal space becomes $\aleph_1$-scwH.
\end{thm}

We may now prove Theorem \ref{thm:FCLC}.
\proof[Proof of Theorem \ref{thm:FCLC}]
   As before, our assumptions imply that $X$ is first countable.
   Assume that $X$ is non-compact.
   Since $X$ is linearly H-closed and normal, $X$ is countably compact.
   By Lemma \ref{lemma:loccomp}, we may assume that $X = \wb{U}$ where $U$ is open and $\sigma$-compact.
   {\bf CW} (Theorem \ref{thm:N-CWH}) shows that $X$ is hereditarily $\aleph_1$-scwH after forcing with a Suslin tree, hence in particular in any model of 
   {\bf MA$_{\omega_1}$(S)[S]} or {\bf PFA(S)[S]}. Since $U$ is open Lindel\"of, 
   Lemma \ref{lemma:Nyikos+} shows that $\wb{U}-U = X - U$ has countable spread. 
   Then Theorem \ref{lemma:MALC} implies 
   that $X-U$ is hereditarily Lindel\"of.
   It follows that $X$ is Lindel\"of and hence compact.
\endproof

We now look at our last result. We will show that our assumptions imply that $X$ is locally compact and apply Theorem \ref{thm:FCLC}.
We need the following consequences of {\bf PFA(S)[S]}.
\begin{thm}[{\cite[Thm 3.5]{DowTall:2016}, \cite[Thm 2.1]{AbrahamTodorcevic:1984}, \cite[Lemma 11]{Tall:PFAArhan}}]
   \label{thm:PFAseparable}
   \
   \\
   {\bf (PFA(S)[S])}
   Separable, hereditarily normal, countably compact spaces are compact.
\end{thm}
(We note in passing that under {\bf PFA}, the theorem holds without `hereditarily', this is 
due to Balogh, Dow, Fremlin and Nyikos \cite[Cor. 2]{Balogh-Dow-Fremlin-Nyikos:1988} and is
the main ingredient in the proof of Fact \ref{fact:PFA}.)
\begin{lemma}[{\bf PFA(S)[S]}]\label{lemma:Nyikos++}
   Let $X$ be first countable, linearly H-closed and hereditarily normal.
   Let $U$ be open and Lindel\"of.
   Then $\wb{U}$ is compact and $\wb{U}-U$ is hereditarily separable.
\end{lemma}
\proof
   $X$ and any closed subset of $X$ are countably compact.
   By {\bf CW} (Theorem \ref{thm:N-CWH}), $X$ is hereditarily $\aleph_1$-scwH.
   By Lemma \ref{lemma:Nyikos+}, $\partial U = \wb{U}-U$ has countable spread.
   By Lemmas \ref{lemma:MAorMA}--\ref{lemma:spreadL}, $\partial U$ is (hereditarily) separable.
   By Theorem \ref{thm:PFAseparable}, $\partial U$ is compact and hence $\wb{U} = U \cup \partial U$
   is Lindel\"of and thus compact.
\endproof

The proof of Theorem \ref{thm:PFASS} is now a formality.
\proof[Proof of Theorem \ref{thm:PFASS}]
  The space is first countable, hence
  by Lemma \ref{lemma:Nyikos++} it is locally compact and Theorem \ref{thm:FCLC} does the rest of the job.
\endproof

This finishes the proofs, the section and the note.


{\footnotesize
\vskip .4cm
\noindent
Mathieu Baillif \\
Haute \'ecole du paysage, d'ing\'enierie et d'architecture (HEPIA) \\
Ing\'enierie des technologies de l'information TIC\\
Gen\`eve -- Suisse
}

\end{document}